\documentclass[12pt,reqno]{amsart}
\usepackage{amsmath,amssymb,latexsym,esint,cite,mathrsfs}
\usepackage{verbatim,wasysym}
\usepackage[left=2.6cm,right=2.6cm,top=3cm,bottom=3cm]{geometry}
\makeatletter \@addtoreset{equation}{section} \makeatother

\usepackage{graphicx}
\usepackage{subfigure}
\usepackage{graphicx,epic,eepic}

\usepackage{color,enumitem,graphicx}
\usepackage[colorlinks=true,urlcolor=blue,
citecolor=red,linkcolor=blue,linktocpage,pdfpagelabels,
bookmarksnumbered,bookmarksopen]{hyperref}

\numberwithin{equation}{section}
\newtheorem{theorem}{Theorem}[section]
\newtheorem{lemma}[theorem]{Lemma}

\newtheorem{remark}[theorem]{Remark}
\newtheorem{corollary}[theorem]{Corollary}
\numberwithin{equation}{section}

\allowdisplaybreaks

\begin{document}

\title[$L^p$-Sobolev inequality]
{A note on the $L^p$-Sobolev inequality}

\author[S. Deng]{Shengbing Deng}
\address{\noindent Shengbing Deng
\newline
School of Mathematics and Statistics, Southwest University,
Chongqing 400715, People's Republic of China}\email{shbdeng@swu.edu.cn}

\author[X. Tian]{Xingliang Tian$^{\ast}$}
\address{\noindent Xingliang Tian  \newline
School of Mathematics and Statistics, Southwest University,
Chongqing 400715, People's Republic of China.}\email{xltian@email.swu.edu.cn}

\thanks{$^{\ast}$ Corresponding author}



\allowdisplaybreaks

\begin{abstract}
The usual Sobolev inequality in $\mathbb{R}^N$, asserts that $\|\nabla u\|_{L^p(\mathbb{R}^N)}
\geq \mathcal{S}\|u\|_{L^{p^*}(\mathbb{R}^N)}$ for $1<p<N$ and $p^*=\frac{pN}{N-p}$, with $\mathcal{S}$ being the sharp constant.
Based on a recent work of Figalli and Zhang [Duke Math. J., 2022], a weak norm remainder term of Sobolev inequality in a subdomain $\Omega\subset \mathbb{R}^N$ with finite measure is established, i.e., for $\frac{2N}{N+1}<p<N$ there exists a constant $\mathcal{C}>0$ independent of $\Omega$ such that
\[
\|\nabla u\|^p_{L^p(\Omega)}
    -\mathcal{S}^p\|u\|^p_{L^{p^*}(\Omega)}
\geq \mathcal{C}|\Omega|^{-\frac{\gamma}{p^*(p-1)}}
\|u\|_{L^{\bar{p}}_w(\Omega)}^{\gamma}\| u\|_{L^{p^*}(\Omega)}^{p-\gamma},\quad \mbox{for all}\ u\in C^\infty_0(\Omega)\setminus\{0\},
\]
where $\gamma=\max\{2,p\}$, $\bar{p}=p^*(p-1)/p$, and $\|\cdot\|_{L^{\bar{p}}_w(\Omega)}$ denotes the weak $L^{\bar{p}}$-norm. Moreover, we establish a sharp upper bound of Sobolev inequality in $\mathbb{R}^N$.
\end{abstract}

\vspace{3mm}

\maketitle

\section{{\bfseries Introduction}}\label{sectir}

    Given $N\geq 2$ and $p\in (1,N)$, denote the homogeneous Sobolev space $\mathcal{D}^{1,p}_0(\mathbb{R}^N)$ be the closure of $C^\infty_c(\mathbb{R}^N)$ with respect to the norm $\|\nabla u\|_{L^p(\mathbb{R}^N)}
    =\left(\int_{\mathbb{R}^N}|\nabla u|^p\mathrm{d}x\right)^{1/p}$.
    The Sobolev inequality states as
    \begin{equation}\label{bzsi}
    \|\nabla u\|_{L^p(\mathbb{R}^N)}\geq \mathcal{S}\|u\|_{L^{p^*}(\mathbb{R}^N)},\quad \mbox{for all}\ u\in \mathcal{D}^{1,p}_0(\mathbb{R}^N),
    \end{equation}
    with $\mathcal{S}=\mathcal{S}(N,p)>0$ being the sharp constant, where $p^*:=\frac{pN}{N-p}$.
    It is well known that
Aubin \cite{Au76} and Talenti \cite{Ta76} found the optimal constant and the extremal functions for \eqref{bzsi}. Indeed, equality is achieved precisely by the functions $cU_{\lambda,z}(x)=c\lambda^{\frac{N-p}{p}}U(\lambda(x-z))$ for all $c\in\mathbb{R}$, $\lambda>0$ and $z\in\mathbb{R}^N$, where
    \[
    U(x)=\gamma_{N,p}(1+|x|^{\frac{p}{p-1}})^{-\frac{N-p}{p}},\quad \mbox{for some constant}\ \gamma_{N,p}>0,
    \]
    which solve the related Sobolev critical equation
    \begin{equation}\label{Ppwh}
    -\mathrm{div}(|\nabla u|^{p-2}\nabla u)= u^{p^*-1} ,\quad u>0 \quad \mbox{in}\ \mathbb{R}^N,\quad u\in \mathcal{D}^{1,p}_0(\mathbb{R}^N),
    \end{equation}
    see \cite{Sc16} for details. Define the set of extremal functions as
    \[
    \mathcal{M}:=\{cU_{\lambda,z}: c\in\mathbb{R}, \lambda>0, z\in\mathbb{R}^N\}.
    \]

    For each bounded domain $\Omega\subset \mathbb{R}^N$, let us define
    \[
    \mathcal{S}(\Omega):=\inf_{u\in \mathcal{D}^{1,p}_0(\Omega)\setminus\{0\}}\frac{\|\nabla u\|_{L^p(\Omega)}}{\|u\|_{L^{p^*}(\Omega)}}.
    \]
    It is well known that $\mathcal{S}(\Omega)=\mathcal{S}(\mathbb{R}^N)=\mathcal{S}$, and $\mathcal{S}(\Omega)$ is never achieved then it is natural to consider the remainder terms. For $p=2$, Br\'{e}zis and Nirenberg \cite{BN83} proved that if $s<\frac{N}{N-2}$ then there is $A=A(\Omega,N,s)>0$ such that
    \begin{equation}\label{rtbn}
    \|\nabla u\|^2_{L^2(\Omega)}
    -\mathcal{S}^2\|u\|^2_{L^{2^*}(\Omega)}
    \geq A\|u\|^2_{L^{s}(\Omega)},\quad \mbox{for all}\ u\in \mathcal{D}^{1,2}_0(\Omega).
    \end{equation}
    Furthermore, the result is sharp in the sense that it is not true if $s=\frac{N}{N-2}$. However, the following refinement is proved by Br\'{e}zis and Lieb \cite{BL85} that
    \begin{equation}\label{rtbl}
    \|\nabla u\|^2_{L^2(\Omega)}
    -\mathcal{S}^2\|u\|^2_{L^{2^*}(\Omega)}
    \geq A'\|u\|^2_{L^{\frac{N}{N-2}}_w(\Omega)},\quad \mbox{for all}\ u\in \mathcal{D}^{1,2}_0(\Omega),
    \end{equation}
    where $\|\cdot\|_{L^s_w(\Omega)}$ denotes the weak $L^s$-norm as
    \begin{align}\label{defwn}
    \|\cdot\|_{L^s_w(\Omega)}:=
    \sup\limits_{D\subset\Omega, |D|>0}|D|^{-\frac{s-1}{s}}
    \int_{D}|\cdot|\mathrm{d}x.
    \end{align}
    Here $|D|$ denotes the Lebesgue measure of $D$. Note that this weak $L^s$-norm is equivalent to the classical weak $L^s$-norm for $s>1$, i.e.,
    \[
    u\in L^s_w(\Omega)\quad \mbox{if and only if}\quad \sup_{t>0} t\mu\{x
    \in\Omega: |u(x)|>t\}^{1/s}<\infty,
    \]
    furthermore, for any $0<t<s$ and $s>1$ with $u\in L^s_w(\Omega)$, we have $\|u\|_{L^t(\Omega)}\leq C_{t,s}\|u\|_{L^s_w(\Omega)}$ which implies the result of \eqref{rtbl} is stronger than \eqref{rtbn}, see \cite[Chapter 5]{CR16} for details.
    Br\'{e}zis and Lieb \cite{BL85} asked a famous question whether a remainder term -- proportional to the quadratic distance of the function $u$ to be the optimizers manifold  $\mathcal{M}$ -- can be added to the right hand side of (\ref{bzsi}). This question was answered affirmatively by Bianchi and Egnell \cite{BE91} by using spectral estimate combined with Lions' concentration and compactness theorem (see \cite{Li85-1}), which reads that there is $c_{\mathrm{BE}}>0$ such that
    \begin{equation}\label{rtbe}
    \|\nabla u\|^2_{L^2(\mathbb{R}^N)} -\mathcal{S}^2\|u\|^2_{L^{2^*}(\mathbb{R}^N)}
    \geq c_{\mathrm{BE}} \inf_{v\in \mathcal{M}}\|\nabla (u-v)\|^2_{L^2(\mathbb{R}^N)}, \quad \mbox{for all}\ u\in \mathcal{D}^{1,2}_0(\mathbb{R}^N),
    \end{equation}
    which can be regarded as a quantitative form of Lion's theorem. Besides, based on the result \eqref{rtbe}, Bianchi and Egnell \cite{BE91} gave a simpler proof of \eqref{rtbl} by showing
    \[
    \|u\|_{L^{\frac{N}{N-2}}_w(\Omega)}\leq C \inf_{v\in \mathcal{M}}\|\nabla (u-v)\|_{L^2(\mathbb{R}^N)}.
    \]
    Chen, Frank and Weth \cite{CFW13} extended \eqref{rtbe} into fractional-order and established \eqref{rtbl} type inequality in a general subdomain $\Omega\subset \mathbb{R}^N$ with $|\Omega|<\infty$. For the general $p\in (1,N)$, Egnell et al. \cite{EPT89} obtained a result of \eqref{rtbn} type that
    \begin{equation}\label{rtbpt}
    \|\nabla u\|^p_{L^p(\Omega)}
    -\mathcal{S}^p\|u\|^p_{L^{p^*}(\Omega)}
    \geq A\|u\|^p_{L^s(\Omega)},\quad \mbox{for all}\ u\in \mathcal{D}^{1,p}_0(\Omega),
    \end{equation}
    for each $s<\bar{p}:=p^*(p-1)/p$, furthermore, the inequality fails if $s=\bar{p}$. For this reason, the number $\bar{p}$ is usually called the critical remainder exponent. Furthermore, Bianchi and Egnell \cite{BE91} conjectured that for all $1<p<N$,
    \begin{align}\label{bec}
    \|\nabla u\|^p_{L^p(\Omega)}
    -\mathcal{S}^p\|u\|^p_{L^{p^*}(\Omega)}
    \geq \mathcal{C}\|u\|^p_{L^{\bar{p}}_w(\Omega)},\quad \mbox{for all}\ u\in \mathcal{D}^{1,p}_0(\Omega),
    \end{align}
    for some $\mathcal{C}>0$. Note that if $1<p\leq \frac{2N}{N+1}$, then $\bar{p}\leq 1$, thus from the definition of weak norm \eqref{defwn} we have $\|u\|_{L^{\bar{p}}_w(\Omega)}
    =|\Omega|^{\frac{1-\bar{p}}{\bar{p}}}
    \|u\|_{L^{1}(\Omega)}$, and the weak norm makes no sense. Therefore, combining with \eqref{rtbpt} we know \eqref{bec} may holds only if $\frac{2N}{N+1}<p<N$.

    When the domain is chosen to be the whole space $\mathbb{R}^N$, Cianchi et al. \cite{CFMP09} first proved a stability version of Lebesgue-type for all $1<p<N$, Figalli and Neumayer \cite{FN19} proved the gradient stability for the Sobolev inequality when $p\geq 2$, Neumayer \cite{Ne20} extended the result in \cite{FN19} to all $1<p<N$. Recently, Figalli and Zhang \cite{FZ22} obtained the sharp stability of Sobolev inequality \eqref{bzsi} for all $1<p<N$, i.e., there is $c_{\mathrm{FZ}}>0$ such that
    \begin{equation}\label{defcbepc}
    \frac{\|\nabla u\|_{L^p(\mathbb{R}^N)}}{\|u\|_{L^{p^*}(\mathbb{R}^N)}}
    -\mathcal{S}
    \geq c_{\mathrm{FZ}} \inf_{v\in \mathcal{M}}\left(\frac{\|\nabla (u-v)\|_{L^p(\mathbb{R}^N)}}{\|\nabla u\|_{L^p(\mathbb{R}^N)}}\right)^\gamma,\quad \mbox{for all}\  u\in \mathcal{D}^{1,p}_0(\mathbb{R}^N)\setminus\{0\},
    \end{equation}
    furthermore, the exponent $\gamma:=\max\{2,p\}$ is sharp. In fact, Figalli and Zhang proved the following equivalent form
    \begin{equation}\label{defcbepcb}
    \|\nabla u\|^p_{L^p(\mathbb{R}^N)}
    -\mathcal{S}^p\|u\|^p_{L^{p^*}(\mathbb{R}^N)}
    \geq c'_{\mathrm{FZ}} \inf_{v\in \mathcal{M}}\|\nabla (u-v)\|_{L^p(\mathbb{R}^N)}^{\gamma}\|\nabla u\|^{p-\gamma}_{L^p(\mathbb{R}^N)}.
    \end{equation}
    When $1<p<2$, \eqref{defcbepcb} looks like a degenerate stability result as in \cite{FP24}.

    As mentioned above, it is natural to consider the weak norm remainder term of $L^p$-Sobolev inequality of \eqref{bec} type which is mentioned by Bianchi and Egnell \cite{BE91}. Recently, Zhou and Zou in \cite[Corollary 1.8]{ZZ23} established the remainder term inequality with weak norm when $\sqrt{N}\leq p<N$, under some assumptions on domain.
    In present paper, based on the sharp stability result \eqref{defcbepc} and the arguments as those in \cite{CFW13}, we consider it in a general subdomain $\Omega\subset \mathbb{R}^N$ with continuous boundary satisfying $|\Omega|<\infty$.

    \begin{theorem}\label{thmprtp}
    Assume $N\geq 2$, $\frac{2N}{N+1}<p<N$, and let $\Omega\subset \mathbb{R}^N$ with continuous boundary satisfy $|\Omega|<\infty$. There exists a constant $\mathcal{C}=\mathcal{C}(N,p)>0$ independent of $\Omega$ such that for all $u\in \mathcal{D}^{1,p}_0(\Omega)\setminus\{0\}$,
    \begin{align}\label{hsrt}
    \|\nabla u\|^p_{L^p(\Omega)}
    -\mathcal{S}^p\|u\|^p_{L^{p^*}(\Omega)}
    \geq \mathcal{C}|\Omega|^{-\frac{\gamma}{p^*(p-1)}}
    \|u\|_{L^{\bar{p}}_w(\Omega)}^{\gamma}
    \|u\|^{p-\gamma}_{L^{p^*}(\Omega)},
    \end{align}
    where $\gamma:=\max\{2,p\}$, $\bar{p}=p^*(p-1)/p$, and $\|\cdot\|_{L^{\bar{p}}_w(\Omega)}$ denotes the weak $L^{\bar{p}}$-norm as in \eqref{defwn}.
    \end{theorem}

    \begin{remark}\label{remp}\rm
    Note that the condition $\frac{2N}{N+1}<p<N$ indicates $\bar{p}=p^*(p-1)/p>1$, then we have $U\in L^{\bar{p}}_w(\mathbb{R}^N)$ (this can be easily verified) which is crucial for comparing $\|u\|_{L^{\bar{p}}_w(\Omega)}$ with $\inf\limits_{v\in \mathcal{M}}\left\|\nabla (u-v)
    \right\|_{L^p(\mathbb{R}^N)}$ (see \eqref{rtbdlbcdl}), however, $\|U\|_{L^{\bar{p}}_w(\mathbb{R}^N)}=+\infty$ if $1<p\leq \frac{2N}{N+1}$. Note also that our result \eqref{hsrt} holds for all $\frac{2N}{N+1}<p<N$, and $\frac{2N}{N+1}<\sqrt{N}$ which indicates our region for $p$ is slightly better than Zhou and Zou
    \cite[Corollary 1.8]{ZZ23}.
    \end{remark}

    From \cite[Theorem 5.16 (a)]{CR16} we know that for any $0<t<s$ and $s>1$,
    \[
    \|u\|_{L^t(\Omega)}\leq \left(\frac{s}{s-t}\right)^{1/t}
    |\Omega|^{\frac{s-t}{st}}\|u\|_{L^s_w(\Omega)}, \quad \mbox{for all}\ u\in L^s_w(\Omega).
    \]
    Then as a direct corollary of Theorem \ref{thmprtp}, we obtain the following Br\'{e}zis and Nirenberg type inequality which can be regarded another form of \eqref{rtbpt}:

    \begin{corollary}\label{thmrtp}
    Assume $N\geq 2$, $\frac{2N}{N+1}<p<N$, and let $\Omega\subset \mathbb{R}^N$ with continuous boundary satisfy $|\Omega|<\infty$. Then for each $t\in (0, \bar{p})$ with $\bar{p}=p^*(p-1)/p$, there exists a constant $\mathcal{C}'=\mathcal{C}'(N,p,t)>0$ independent of $\Omega$ such that for all $u\in \mathcal{D}^{1,p}_0(\Omega)\setminus\{0\}$,
    \begin{align}\label{hsrts}
    \|\nabla u\|^p_{L^p(\Omega)}
    -\mathcal{S}^p\|u\|^p_{L^{p^*}(\Omega)}
    \geq \mathcal{C}'
    |\Omega|^{-\frac{\gamma(p^*-t)}{tp^*}}
    \|u\|_{L^{t}(\Omega)}^{\gamma}
    \|u\|^{p-\gamma}_{L^{p^*}(\Omega)},
    \end{align}
    where $\gamma:=\max\{2,p\}$.
    \end{corollary}

    Finally, following the arguments as those in the recent work \cite{DTYZ24}, we give a upper bound of Sobolev inequality in $\mathbb{R}^N$, which may has its own interests.

    \begin{theorem}\label{thmprtpl}
    Assume $1<p<N$. There exists a constant  $\mathcal{C}''=\mathcal{C}''(N,p)>0$ such that for all $u\in \mathcal{D}^{1,p}_0(\mathbb{R}^N)\setminus\{0\}$,
    \begin{align}\label{defcbepcb2}
    \|\nabla u\|^p_{L^p(\mathbb{R}^N)}
    -\mathcal{S}^p\|u\|^p_{L^{p^*}(\mathbb{R}^N)}
    \leq \mathcal{C}''\inf_{v\in \mathcal{M}}\|\nabla (u-v)\|_{L^p(\mathbb{R}^N)}^{\zeta}\|\nabla u\|^{p-\zeta}_{L^p(\mathbb{R}^N)},
    \end{align}
    furthermore, the exponent $\zeta:=\min\{2,p\}$ is sharp.
    \end{theorem}

    \begin{remark}\label{remupsies}\rm
    The sharpness of the exponent $\zeta=\min\{2,p\}$ in \eqref{defcbepcb2} follows directly from \cite[Remark 1.2]{FZ22}.
    \end{remark}

    The paper is organized as follows: in Section \ref{sectpmr}, we give the proof of weak norm remainder term of Sobolev inequality in a general subdomain $\Omega\subset \mathbb{R}^N$ with $|\Omega|<\infty$. Section \ref{sectubsi} is devoted to proving the upper bound of Sobolev inequality in whole space $\mathbb{R}^N$.

\section{Sobolev inequality with remainder terms in a subdomain}\label{sectpmr}

In order to prove \eqref{hsrt}, by homogeneity we can always assume that $\|u\|_{L^{p^*}(\Omega)}=1$.
Note that $|\nabla u|\geq |\nabla |u||$ thus it is suffices to consider $|u|$ instead of $u$ in \eqref{hsrt}.
    By the rearrangement inequality, we have
    \begin{align*}
    \|\nabla u^*\|_{L^p(B_R)}\leq \|\nabla u\|_{L^p(\Omega)},\quad \|u^*\|_{L^{p^*}(B_R)}=\|u\|_{L^{p^*}(\Omega)}, \quad
    \|u^*\|_{L^{\bar{p}}_w(B_R)}= \|u\|_{L^{\bar{p}}_w(\Omega)},
    \end{align*}
    where $\|\cdot\|_{L^{\bar{p}}_w(\Omega)}$ denotes the weak $L^{\bar{p}}$-norm as in \eqref{defwn} with $\bar{p}=p^*(p-1)/p$, and $u^*$ denotes the symmetric decreasing rearrangement of nonnegative function $u$ extended to zero outside $\Omega$, and
    \[
    |\Omega|=|B_R|\quad \mbox{for some}\ R\in (0,\infty),\quad B_R:=B(\mathbf{0},R).
    \]
    Moreover, by using H\"{o}lder inequality we have
    \[
    \|u\|_{L^{\bar{p}}_w(\Omega)}\leq \|u\|_{L^{\bar{p}}(\Omega)}\leq \|u\|_{L^{p^*}(\Omega)}|\Omega|^{\frac{1}{p\cdot\bar{p}}}
    =|\Omega|^{\frac{1}{p^*(p-1)}}.
    \]
Therefore, it is sufficient to consider the case in which $\Omega$ is a ball of radius $R$ at origin and $u$ is nonnegative symmetric decreasing, i.e.,
    \begin{equation}\label{rtbdlb}
    \|\nabla u\|^p_{L^p(\Omega)}
    -\mathcal{S}^p\|u\|^p_{L^{p^*}(\Omega)}
    \geq C |B_R|^{-\frac{\gamma}{p^*(p-1)}}
    \|u\|_{L^{\bar{p}}_w(B_R)}^{\gamma},
    \end{equation}
    for all $u\in \mathfrak{R}_0^{1,p}(B_R)$ satisfying
    \begin{align}\label{rtzhc}
    \|u\|_{L^{p^*}(B_R)}=1,\quad \|\nabla u\|^p_{L^p(\Omega)}
    -\mathcal{S}^p\|u\|^p_{L^{p^*}(\Omega)}\ll 1,
    \end{align}
    where $\gamma=\max\{p,2\}$, and $\mathfrak{R}_0^{1,p}(B_R)$ consists all nonnegative and radial functions in $\mathcal{D}^{1,p}_0(B_R)$ with support in closed ball $\overline{B_R}$.
    Note that \eqref{rtzhc} implies $\|\nabla u\|_{L^p(B_R)}$ is bounded away from zero and infinity, i.e., $c_0\leq \|\nabla u\|_{L^p(B_R)}\leq C_0$ for some constants $C_0\geq c_0>0$. Therefore, \eqref{rtbdlb} is equivalent to
    \begin{equation}\label{rtbdlbn}
    \|\nabla u\|^p_{L^p(\Omega)}
    -\mathcal{S}^p\|u\|^p_{L^{p^*}(\Omega)}
    \geq C'|B_R|^{-\frac{\gamma}{p^*(p-1)}}
    \|u\|_{L^{\bar{p}}_w(B_R)}^{\gamma}
    \|\nabla u\|^{p-\gamma}_{L^{p}(\Omega)},
    \end{equation}
    for all $u\in \mathfrak{R}_0^{1,p}(B_R)$ satisfying \eqref{rtzhc}. Then, the remainder inequality \eqref{rtbdlbn} will follow immediately from the following lemma and \eqref{defcbepcb}.

    \begin{lemma}\label{proprtbd}
    Assume $N\geq 2$ and $\frac{2N}{N+1}<p<N$. There exists a constant $\mathcal{B}>0$ depending only on $N$ and $p$ such that for all $u\in \mathfrak{R}_0^{1,p}(B_R)$ satisfying \eqref{rtzhc},
    \begin{equation}\label{rtbdke}
    \|u\|_{L^{\bar{p}}_w(B_R)}
    \leq \mathcal{B} |B_R|^{\frac{1}{p^*(p-1)}}\inf_{v\in \mathcal{M}}\|\nabla u-\nabla v\|_{L^p(\mathbb{R}^N)}.
    \end{equation}
    \end{lemma}

    \begin{proof}
    We follow the arguments as those in \cite[Proposition 3]{CFW13}. Let $u\in \mathfrak{R}_0^{1,p}(B_R)$ satisfy \eqref{rtzhc}. Firstly, we notice that \eqref{defcbepcb} and \eqref{rtzhc} indicate
    \[
    \inf_{v\in \mathcal{M}}\|\nabla u-\nabla v\|_{L^p(\mathbb{R}^N)}\ll 1,
    \]
    then from
    \cite[Lemma 4.1]{FZ22} we know that $\inf\limits_{v\in \mathcal{M}}\|\nabla u-\nabla v\|_{L^p(\mathbb{R}^N)}$ can always be attained, i.e.,
    \[
    \inf_{v\in \mathcal{M}}\|\nabla u-\nabla v\|_{L^p(\mathbb{R}^N)}=\|\nabla (u-cU_{\lambda,0})\|_{L^p(\mathbb{R}^N)}\quad  \mbox{for some}\ c\in\mathbb{R}, \ \lambda>0,
    \]
    thanks to $u$ is radially symmetric.
    Furthermore, since $u$ is nonnegative, we have $c>0$.

    As stated previous, \eqref{rtzhc} implies $\|\nabla u\|_{L^p(B_R)}$ is bounded away from zero and infinity, i.e., $c_0\leq \|\nabla u\|_{L^p(B_R)}\leq C_0$ for some constants $C_0\geq c_0>0$. Let $\rho\in (0,c_0)$ be given by
    \begin{align}\label{defrho}
    \frac{\rho \|\nabla U\|_{L^p(\mathbb{R}^N)}}{(c_0-\rho)\mathcal{S}}
    =\gamma_{N,p}\left(|\mathbb{S}^{N-1}|\int^\infty_1 \frac{r^{N-1}}{(1+r^{\frac{p}{p-1}})
    ^{N}}\mathrm{d}r
    \right)^{1/p^*},
    \end{align}
    where $\gamma_{N,p}=U(0)$.
    So
    \[
    \inf_{v\in \mathcal{M}}\|\nabla u-\nabla v\|_{L^p(\mathbb{R}^N)}<\rho,
    \]
    due to $\inf\limits_{v\in \mathcal{M}}\|\nabla u-\nabla v\|_{L^p(\mathbb{R}^N)}\ll 1$ and $\rho\in (0,c_0)$ is a fixed constant. Note that
    \begin{align*}
    |\|\nabla u\|_{L^p(\mathbb{R}^N)}-\|\nabla (cU_{\lambda,0})\|_{L^p(\mathbb{R}^N)}|
    \leq \|\nabla (u-cU_{\lambda,0})\|_{L^p(\mathbb{R}^N)}
    =\inf_{v\in \mathcal{M}}\|\nabla u-\nabla v\|_{L^p(\mathbb{R}^N)}<\rho,
    \end{align*}
    which implies
    \[
    \frac{c_0-\rho}{\|\nabla U\|_{L^p(\mathbb{R}^N)}}
    \leq \frac{\|\nabla u\|_{L^p(\mathbb{R}^N)}-\rho}{\|\nabla U\|_{L^p(\mathbb{R}^N)}}\leq c\leq \frac{\|\nabla u\|_{L^p(\mathbb{R}^N)}+\rho}{\|\nabla U\|_{L^p(\mathbb{R}^N)}}
    \leq \frac{C_0+\rho}{\|\nabla U\|_{L^p(\mathbb{R}^N)}}.
    \]
    Then we have
    \begin{align}\label{dug}
    \inf_{v\in \mathcal{M}}\|\nabla u-\nabla v\|_{L^p(\mathbb{R}^N)}
    & =\|\nabla (u-cU_{\lambda,0})\|_{L^p(\mathbb{R}^N)}
    \nonumber\\
    & \geq \mathcal{S} \|u-cU_{\lambda,0}\|_{L^{p^*}
    (\mathbb{R}^N)}
    \nonumber\\
    & \geq \mathcal{S} c\|U_{\lambda,0}\|_{L^{p^*}(\mathbb{R}^N\setminus B_R)}
    \nonumber\\
    & \geq \left(\frac{c_0-\rho}
    {\|\nabla U\|_{L^p(\mathbb{R}^N)}}\right)\mathcal{S} \|U_{\lambda,0}\|_{L^{p^*}(\mathbb{R}^N\setminus B_R)},
    \end{align}
    hence
    \begin{align}\label{Uwgo}
    \|U_{\lambda,0}\|^{p^*}_{L^{p^*}(\mathbb{R}^N\setminus B_R)}
    & \leq \left(\frac{\inf\limits_{v\in \mathcal{M}}\|\nabla u-\nabla v\|_{L^p(\mathbb{R}^N)}\|\nabla U\|_{L^p(\mathbb{R}^N)}}
    {(c_0-\rho)\mathcal{S}}\right)^{p^*}
    \nonumber\\
    & \leq \left(\frac{\rho \|\nabla U\|_{L^p(\mathbb{R}^N)}}{(c_0-\rho)
    \mathcal{S}}\right)^{p^*}
    =\gamma_{N,p}^{p^*}|\mathbb{S}^{N-1}|\int^\infty_1 \frac{r^{N-1}}{(1+r^{\frac{p}{p-1}})
    ^N}\mathrm{d}r
    \end{align}
    by the choice of $\rho$ in \eqref{defrho}. On the other hand, we compute
    \begin{align*}
    \|U_{\lambda,0}\|^{p^*}_{L^{p^*}(\mathbb{R}^N\setminus B_R)}
    & =\gamma_{N,p}^{p^*}|\mathbb{S}^{N-1}|\int^\infty_R \frac{r^{N-1}\lambda^{N}}{(1+(\lambda r)^{\frac{p-\beta}{p-1}})
    ^N}\mathrm{d}r
    \\
    & =\gamma_{N,p}^{p^*}|\mathbb{S}^{N-1}|\int^\infty_{\lambda R} \frac{r^{N-1}}{(1+r^{\frac{p}{p-1}})
    ^N}\mathrm{d}r,
    \end{align*}
    which implies $\lambda R\geq 1$ and therefore
    \begin{align}\label{Uwg}
    \|U_{\lambda,0}\|^{p^*}_{L^{p^*}(\mathbb{R}^N\setminus B_R)}
    & \geq 2^{-N}\gamma_{N,p}^{p^*}
    |\mathbb{S}^{N-1}|\int^\infty_{\lambda R} r^{-\frac{N}{p-1}-1}\mathrm{d}r
    \nonumber\\
    & =2^{-N}
    \gamma_{N,p}^{p^*}|\mathbb{S}^{N-1}|\frac{N}{p-1}
    R^{-\frac{N}{p-1}}
    \lambda^{-\frac{N}{p-1}}.
    \end{align}
    Combining \eqref{Uwgo} and \eqref{Uwg}, from \eqref{dug}, we conclude that
    \begin{align}\label{dugf}
    \inf_{v\in \mathcal{M}}\|\nabla u-\nabla v\|_{L^p(\mathbb{R}^N)}\geq \underline{C}R^{-\frac{N-p}{p(p-1)}}\lambda^{-\frac{N-p}{p(p-1)}}
    \end{align}
    with
    \[
    \underline{C}:=\frac{(c_0-\rho)\mathcal{S}
    \gamma_{N,p}}
    {\|\nabla U\|_{L^p(\mathbb{R}^N)}}
    \left(2^{-N}
    |\mathbb{S}^{N-1}|
    \frac{N}{p-1}\right)^{1/p^*},
    \]
    thanks to $\|\nabla U\|_{L^p(\mathbb{R}^N)}= \mathcal{S}\|U\|_{L^{p^*}(\mathbb{R}^N)}$ and $\|\nabla U\|^p_{L^p(\mathbb{R}^N)}=\|U\|^{p^*}_{L^{p^*}(\mathbb{R}^N)}$ imply $\|\nabla U\|_{L^p(\mathbb{R}^N)}=\mathcal{S}^{\frac{p^*}{p^*-p}}$.
    Then we have
    \begin{align}\label{rtbdlbcdl}
    \|u\|_{L^{\bar{p}}_w(B_R)}
    & \leq \|cU_{\lambda,0}\|_{L^{\bar{p}}_w(B_R)}
    + \|u-cU_{\lambda,0}\|_{L^{\bar{p}}_w(B_R)}
    \nonumber \\
    & \leq c\lambda^{-\frac{N-p}{p(p-1)}}
    \|U\|_{L^{\bar{p}}_w(B_{\lambda R})}
    + \|u-cU_{\lambda,0}\|_{L^{\bar{p}}(B_R)}
    \nonumber \\
    & \leq \frac{C_0+\rho}{\|\nabla U\|_{L^p(\mathbb{R}^N)}}\lambda^{-\frac{N-p}{p(p-1)}}
    \|U\|_{L^{\bar{p}}_w(\mathbb{R}^N)}
    +|B_R|^{\frac{1}{p\cdot \bar{p}}}\mathcal{S}^{-1}
    \left\|\nabla(u-cU_{\lambda,0})
    \right\|_{L^p(\mathbb{R}^N)}
    \nonumber \\
    & \leq \mathcal{B}|B_R|^{\frac{1}{p^*(p-1)}}
    \inf\limits_{v\in \mathcal{M}}\|\nabla u-\nabla v\|_{L^p(\mathbb{R}^N)}
    \end{align}
    with
    \[
    \mathcal{B}:=\frac{(C_0+\rho)\|U\|_{L^{\bar{p}}_w(\mathbb{R}^N)}}
    {\underline{C}|\mathbb{S}^{N-1}|^{\frac{1}{p^*(p-1)}}
    \mathcal{S}^{\frac{p^*}{p^*-p}}}
    +\mathcal{S}^{-1}.
    \]
    Now, the proof of \eqref{rtbdke} is completed.
    \end{proof}

Now, we are ready to prove the weak-Lebesgue remainder inequality \eqref{hsrt}.

\vskip0.25cm

\noindent{\bf \em Proof of Theorem \ref{thmprtp}}.
As stated in the beginning of this section, in order to prove the weak-Lebesgue remainder inequality \eqref{hsrt}, it is sufficient to prove \eqref{rtbdlbn} under the condition \eqref{rtzhc}, which follows directly from Lemma \ref{proprtbd} and \eqref{defcbepcb}.
\qed

\vskip0.25cm

\section{Upper bound of Sobolev inequality in whole space}\label{sectubsi}

    In this section, we consider the upper bound of Sobolev inequality \eqref{bzsi}. In order to do this, firstly, we need the following algebraic inequalities.

    \begin{lemma}\label{lemui1p} {\rm \cite[Lemma A.4]{Sh00}}
    Let $x, y\in\mathbb{R}^N$, the following inequalities hold.
    \begin{itemize}
    \item[$(i)$]
    If $p\geq 2$ then
    \begin{align}\label{uinb1p}
    |x+y|^p
    & \leq |x|^p+ p|x|^{p-2}x\cdot y
    +\frac{p(p-1)}{2}(|x|+|y|)^{p-2}|y|^2.
    \end{align}
    \item[$(ii)$]
    If $1<p<2$ then there exists a constant $\gamma_p>0$ such that
    \begin{align}\label{uinb2p}
    |x+y|^p
    \leq |x|^p+ p|x|^{p-2}x\cdot y
    +\gamma_p|y|^{p}.
    \end{align}
    \end{itemize}
    \end{lemma}

    \begin{lemma}\label{lemui1pr} {\rm \cite[Lemma 2.1]{FZ22}}
    Let $x, y\in\mathbb{R}^N$. Then for any $\kappa>0$, there exists a constant $\mathcal{C}_1=\mathcal{C}_1(r,\kappa)>0$ such that
     the following inequalities hold.

    \begin{itemize}
    \item[$(i)$]
    If $r\geq 2$ then
    \begin{align*}
    |x+y|^r
    & \geq |x|^r+ r|x|^{r-2}x\cdot y+ \frac{1-\kappa}{2}\left(r|x|^{r-2}|y|^2+ r(r-2)|\bar{\omega}|^{r-2}(|x|-|x+y|)^2 \right)
    +\mathcal{C}_1 |y|^r ,
    \end{align*}
    where
    \begin{eqnarray*}
    \bar{\omega}=\bar{\omega}(x,x+y)=
    \left\{ \arraycolsep=1.5pt
       \begin{array}{ll}
        \left(\frac{|x+y|}{|x|}\right)^{\frac{1}{r-2}}(x+y),\ \ &{\rm if}\ \  |x+y|\leq |x|\\[3mm]
        x,\ \ &\mathrm{if}\ |x|<|x+y|
        \end{array}.
    \right.
    \end{eqnarray*}

    \item[$(ii)$] If $1<r<2$ then
    \begin{align*}
    |x+y|^r
    & \geq |x|^r+ r|x|^{r-2}x\cdot y+ \frac{1-\kappa}{2}\left(r|x|^{r-2}|y|^2+ r(r-2)|\tilde{\omega}|^{r-2}(|x|-|x+y|)^2 \right) \nonumber\\
    & \quad +\mathcal{C}_1\min\{|y|^r,|x|^{r-2}|y|^2\},
    \end{align*}
    where
    \begin{eqnarray*}
    \tilde{\omega}=\tilde{\omega}(x,x+y)=
    \left\{ \arraycolsep=1.5pt
       \begin{array}{ll}
        \left(\frac{|x+y|}{(2-r)|x+y|+(r-1)|x|}
        \right)^{\frac{1}{r-2}}x,\ \ &\mathrm{if}\ |x|<|x+y|\\[3mm]
        x,\ \ &\mathrm{if}\  |x+y|\leq |x|
        \end{array}.
    \right.
    \end{eqnarray*}
    \end{itemize}
    \end{lemma}

    Note that if $1<r<2$, then $|x|^{r-2}|y|^2+(r-2)|\tilde{\omega}|^{r-2}(|x|-|x+y|)^2\geq 0$ for any $x\neq 0$, see \cite[(2.2)]{FZ22} for details. Therefore, from Lemma \ref{lemui1pr} we deduce that for each $r>1$,
    \begin{align}\label{uinb2pb}
    |a+b|^r
    & \geq |a|^r+ r|a|^{r-2}ab,\quad \mbox{for all}\ a,b\in\mathbb{R}.
    \end{align}

    The main ingredient of the upper bound of Sobolev inequality is contained in the following lemma, in which the behavior near the extremal functions set $\mathcal{M}$ is studied.

    \begin{lemma}\label{lemma:rtnml}
    Suppose $1<p<N$.
    There exists a constant $\varrho>0$ such that for any sequence $\{u_n\}\subset \mathcal{D}^{1,p}_{0}(\mathbb{R}^N)\backslash \mathcal{M}$ satisfying $\|\nabla u_n\|_{L^p(\mathbb{R}^N)}=1$ and $\inf\limits_{v\in \mathcal{M}}\|\nabla (u_n-v)\|_{L^p(\mathbb{R}^N)}\to 0$,
    \begin{align}\label{rtnmb}
    \limsup_{n\to\infty}
    \frac{1
    -\mathcal{S}^p\|u_n\|^p_{L^{p^*}(\mathbb{R}^N)}}
    {\inf\limits_{v\in \mathcal{M}}\|\nabla (u_n-v)\|_{L^p(\mathbb{R}^N)}^{\zeta}}
    \leq \varrho,
    \end{align}
    where $\zeta=\min\{2,p\}$.
    \end{lemma}

    \begin{proof}
    Since $\|\nabla u_n\|_{L^p(\mathbb{R}^N)}=1$ and $d_n:=\inf\limits_{v\in \mathcal{M}}\|\nabla (u_n-v)\|_{L^p(\mathbb{R}^N)}\to 0$, from
    \cite[Lemma 4.1]{FZ22} we know that $d_n$ can always be attained for each sufficiently large $n$, i.e., there are  $c_n\in\mathbb{R}\setminus\{0\}$, $\lambda_n>0$ and $z_n\in\mathbb{R}^N$ such that $d_n=\|\nabla(u_n-c_nU_{\lambda_n,z_n})\|_{L^p(\mathbb{R}^N)}$. Since $\mathcal{M}$ is a smooth $(N+2)$-manifold and the tangential space at $c_n U_{\lambda_n,z_n}$ is given by
    \begin{align*}
    T_{c_n U_{\lambda_n,z_n}}\mathcal{M}
    =\mathrm{Span}\left\{U_{\lambda_n,z_n},\ \frac{\partial U_{\lambda_n,z_n}}{\partial \lambda_n},\ \frac{\partial U_{\lambda_n,z_n}}{\partial z_n^i},i=1,\ldots,N\right\},
    \end{align*}
    we rewrite $u_n$ as
    \begin{equation}\label{defunwn}
    u_n=c_n U_{\lambda_n,z_n}+d_n w_n,
    \end{equation}
    then $w_n$ is perpendicular to $T_{c_n U_{\lambda_n,z_n}}\mathcal{M}$ satisfying $\|\nabla w_n\|_{L^p(\mathbb{R}^N)}=1$ and
    \begin{align*}
    \int_{\mathbb{R}^N} |\nabla U_{\lambda_n,z_n}|^{p-2}\nabla U_{\lambda_n,z_n}\cdot \nabla w_n \mathrm{d}x
    =\int_{\mathbb{R}^N}U_{\lambda_n,z_n}^{p^*-1}w_n
    \mathrm{d}x
    =0,
    \end{align*}
    thanks to $U_{\lambda_n,z_n}$ is the solution of Sobolev critical equation \eqref{Ppwh}.

    From \eqref{uinb2pb} we have
    \begin{align*}
    \|u_n\|^{p^*}_{L^{p^*}(\mathbb{R}^N)}
    \geq |c_n|^{p^*}\int_{\mathbb{R}^N}|U_{\lambda_n,z_n}|^{p^*} \mathrm{d}x
    +p|c_n|^{p^*-2}c_nd_n \int_{\mathbb{R}^N}U_{\lambda_n,z_n}^{p^*-1}w_n \mathrm{d}x
    = |c_n|^{p^*}\|U\|^{p^*}_{L^{p^*}(\mathbb{R}^N)},
    \end{align*}
    thus
    \begin{align}\label{unp*b}
    \|u_n\|^{p}_{L^{p^*}(\mathbb{R}^N)}
    \geq |c_n|^{p}\|U\|^{p}_{L^{p^*}(\mathbb{R}^N)},
    \end{align}
    When $p\geq 2$, from \eqref{uinb1p} we have
    \begin{align*}
    \|\nabla u_n\|^p_{L^p(\mathbb{R}^N)}
    & = \int_{\mathbb{R}^N}|c_n \nabla U_{\lambda_n,z_n}+d_n \nabla w_n|^p \mathrm{d}x\nonumber\\
    & \leq |c_n|^{p}\int_{\mathbb{R}^N}|\nabla U_{\lambda_n,z_n}|^p \mathrm{d}x
    +p|c_n|^{p-2}c_nd_n \int_{\mathbb{R}^N}|\nabla U_{\lambda_n,z_n}|^{p-2}  \nabla U_{\lambda_n,z_n}\cdot \nabla w_n \mathrm{d}x  \nonumber\\
    & \quad +\frac{p(p-1)}{2}d_n^2\int_{\mathbb{R}^N}\left(|c_n\nabla U_{\lambda_n,z_n}|+ |d_n \nabla w_n|\right)^{p-2}  |\nabla w_n|^2\mathrm{d}x
    \\
    & = |c_n|^{p}\|\nabla U\|^{p}_{L^{p}(\mathbb{R}^N)}
    +\frac{p(p-1)}{2}d_n^2\int_{\mathbb{R}^N}\left(|c_n\nabla U_{\lambda_n,z_n}|+ |d_n \nabla w_n|\right)^{p-2}  |\nabla w_n|^2\mathrm{d}x.
    \end{align*}
    Moreover, for $p\geq 2$, by H\"{o}lder inequality we have
    \begin{align*}
    & \int_{\mathbb{R}^N}\left(|c_n\nabla U_{\lambda_n,z_n}|+ |d_n \nabla w_n|\right)^{p-2}  |\nabla w_n|^2\mathrm{d}x
    \\
    & \leq \left(\int_{\mathbb{R}^N}\left(|c_n\nabla U_{\lambda_n,z_n}|+ |d_n \nabla w_n|\right)^{p}\mathrm{d}x\right)^{\frac{p-2}{p}}
    \left(\int_{\mathbb{R}^N}|\nabla w_n|^p\mathrm{d}x\right)^{\frac{2}{p}}
    \\
    &\leq 2^{\frac{(p-1)(p-2)}{p}}\left(|c_n|^p\int_{\mathbb{R}^N}|\nabla U_{\lambda_n,z_n}|^p\mathrm{d}x
    +d_n^p\int_{\mathbb{R}^N}|\nabla w_n|^p\mathrm{d}x\right)^{\frac{p-2}{p}}
    \\
    & = 2^{\frac{(p-1)(p-2)}{p}}\left(|c_n|^p\|\nabla U\|^{p}_{L^{p}(\mathbb{R}^N)}+d_n^p\right)^{\frac{p-2}{p}},
    \end{align*}
    thanks to $(a+b)^p\leq 2^{p-1}(a^p+b^p)$ for all $a,b\geq 0$ and $p>1$. Since $\|\nabla u_n\|_{L^p(\mathbb{R}^N)}=1$, then from Lemma \ref{lemui1pr} it is not difficult to verify that $|c_n|$ is bounded. Therefore,
    \begin{align}\label{epknugpg2l}
    \|\nabla u_n\|^p_{L^p(\mathbb{R}^N)}
    \leq |c_n|^{p}\|\nabla U\|^{p}_{L^{p}(\mathbb{R}^N)}
    +Cd_n^2.
    \end{align}
    Thus for $p\geq 2$, combing with \eqref{unp*b} and \eqref{epknugpg2l} we have
    \begin{align}\label{epknugpg2lc}
    \|\nabla u_n\|^p_{L^p(\mathbb{R}^N)}
    -\mathcal{S}^p\|u_n\|^{p}_{L^{p^*}(\mathbb{R}^N)}
    \leq |c_n|^{p}\|\nabla U\|^{p}_{L^{p}(\mathbb{R}^N)}+Cd^2
    -|c_n|^{p}\|U\|^{p}_{L^{p^*}(\mathbb{R}^N)}
    =Cd_n^2.
    \end{align}
    When $1<p<2$, from \eqref{uinb2p} we have
    \begin{align}\label{epknugpl2l}
    \|\nabla u_n\|^p_{L^p(\mathbb{R}^N)}
    & \leq |c_n|^{p}\int_{\mathbb{R}^N}|\nabla U_{\lambda_n,z_n}|^p \mathrm{d}x
    +p|c_n|^{p-2}c_nd_n \int_{\mathbb{R}^N}|\nabla U_{\lambda_n,z_n}|^{p-2}  \nabla U_{\lambda_n,z_n}\cdot \nabla w_n \mathrm{d}x  \nonumber\\
    & \quad +\gamma_p d_n^p\int_{\mathbb{R}^N}|\nabla w_n|^p\mathrm{d}x
    \nonumber\\
    & = |c_n|^{p}\|\nabla U\|^{p}_{L^{p}(\mathbb{R}^N)}
    +\gamma_p d_n^p,
    \end{align}
    for some constant $\gamma_p>0$.
    Thus for $1<p<2$, combing with \eqref{unp*b} and \eqref{epknugpl2l} we have
    \begin{align}\label{epknugpl2lc}
    \|\nabla u_n\|^p_{L^p(\mathbb{R}^N)}
    -\mathcal{S}^p\|u_n\|^{p}_{L^{p^*}(\mathbb{R}^N)}
    & \leq |c_n|^{p}\|\nabla U\|^{p}_{L^{p}(\mathbb{R}^N)}
    +\gamma_p d_n^p
    -|c_n|^{p}\mathcal{S}^p\|U\|^{p}_{L^{p^*}(\mathbb{R}^N)}
    \nonumber\\
    & =\gamma_p d_n^p.
    \end{align}
    Therefore, \eqref{rtnmb} follows directly from \eqref{epknugpg2lc} and \eqref{epknugpl2lc}.
    \end{proof}

    Now, we are ready to prove the upper bound of Sobolev inequality.
\vskip0.25cm

\noindent{\bf \em Proof of Theorem \ref{thmprtpl}.}  By homogeneity, we can assume that $\|\nabla u\|_{L^p(\mathbb{R}^N)}=1$. Now, we argue by contradiction. In fact, if the theorem is false then there exists a sequence $\{u_n\}\subset \mathcal{D}^{1,p}_{0}(\mathbb{R}^N)\backslash \mathcal{M}$ satisfying $\|\nabla u_n\|_{L^p(\mathbb{R}^N)}=1$ such that
    \begin{align*}
    \frac{1
    -\mathcal{S}^p\|u_n\|^p_{L^{p^*}(\mathbb{R}^N)}}
    {\inf\limits_{v\in \mathcal{M}}\|\nabla (u_n-v)\|_{L^p(\mathbb{R}^N)}^{\zeta}}
    \to +\infty,\quad \mbox{as}\ n\to \infty,
    \end{align*}
    where $\zeta=\min\{2,p\}$. Since $0\leq 1
    -\mathcal{S}^p\|u_n\|^p_{L^{p^*}(\mathbb{R}^N)}\leq 1$ for $\|\nabla u_n\|_{L^p(\mathbb{R}^N)}=1$, it must be $\inf\limits_{v\in \mathcal{M}}\|\nabla (u_n-v)\|_{L^p(\mathbb{R}^N)}\to 0$ which leads to a contradiction by Lemma \ref{lemma:rtnml}.
\qed

\vskip0.25cm

\noindent{\bfseries Acknowledgements}

The research has been supported by National Natural Science Foundation of China (No. 12371121).

    \end{document}